\documentclass{amsart}
\usepackage{enumerate}
\usepackage{amsmath}
\usepackage{amssymb}
\usepackage{amsthm, amsfonts, amstext}
\usepackage{epsfig, color}
\usepackage{graphicx}
\newtheorem{theorem}{Theorem}
\newtheorem{lemma}{Lemma}
\newtheorem{proposition}{Proposition}
\newtheorem{corollary}{Corollary}
\def\Z{{\mathbb Z}}        % integers
\def\R{{\mathbb R}}        % reals
\def\N{{\mathbb N}}        % natural numbers

\def\X{{\Omega}}            % State space
\def\P{{\mathbb P}}        % probability
\def\E{{\mathbb E}}        % expectation
\def\1{{\mathbf 1}}        % indicator
\def\G{{\mathcal G}}        % some sigma field
\def\F{{\mathcal F}}        % some sigma field

\newcommand\fracc[2]{{#1}/{#2}   }
\newcommand\floor[1]{\lfloor  {#1}   \rfloor}
\newcommand\naiveset[1]{ \{ {#1} \} }
\newcommand\textt[1]{ \  {\text{#1}}  \   }
\newcommand\textte[1]{ \  {\text{#1}}     }
\newcommand\PP[1]{ \P  \{ {#1} \} }

\newcommand\obracket[1] {   \left( {#1} \right) }

\newcommand\PPp[1]{ \P \obracket{#1} }

\def\Ppp{{\mathcal P}}

\newcommand\one[1]{ \mathbf{1}_{ [#1]}    }

\def\e{{ \varepsilon }}

\def\a{ {\alpha     } }
\def\b{ { \beta}   }
\def\Oo{ {\mathcal O     } }

\newcommand{\df}{\bf}

\newcommand\cardsupp[1]{  \# [#1] }
\newcommand\ns[1]{ \left\{ {#1} \right\} }

\def\CC{{\mathcal C}}

\def\clump{{L}}
\def\old{ {2    } }
\def\oldd{ {4    } }
\def\coin{{ \{0,1\}^{\Z^d}}}

\def\choiceofalpha{\frac{1}{1 + \frac{d}{4}}}
\def\choiceofbeta{\frac{2}{1 + \frac{4}{d}}} 
\def\result{ \PPp{Z_{ \Phi} > r}  \leq c (\log r)^{ \oldd  } r^ {-\b}  }
\def\conditions{ \b = \b(d) =   \choiceofbeta   }
\def\resultmm{ \PPp{Z_{ \Phi} > r}  \leq c  r^ {-\frac{1}{2}}   }  
\def\resultalld{ \PPp{Z_{ \Phi} > r}  \leq c  r^ {-\frac{d}{2}}   }

\def\bgset{S(2r_{k})}
\def\cconstantk{c_{1}}
\def\clatticep{c_{2}}
\def\ctimesfirstterm{c_{3}}
\def\cclt{c_{5}}
\def\cklog{c_{4}}
\def\cmethod1{c_{6}}
\def\cestager{c_{7}}
\def\ccutterconst{c_{8}}
\def\constantmethod1{c_{9}}
\def\constantmethoda{c_{_{10}}}
\def\constantmethodb{c_{_{11}}}
\def\constantmethodaa{c_{_{12}}}
\def\constantmethodbb{c_{_{13}}}
\def\ccube{c_{_{14}}}
\def\csidele{c_{_{15}}}
\newcommand{\abs}[1]{ {\Bigg| {#1} \Bigg|}  }

\newcommand{\kseed}{$k$-seed}
\newcommand{\kseeds}{$k$-seeds}
\newcommand{\kblob}{$k$-blob}
\newcommand{\kblobs}{$k$-blobs}

\newcommand{\kcutters}{$k$-cutters}
\newcommand{\kbad}{$k$-bad}

\newcommand{\Kbad}{$K$-bad}

\newcommand{\jseeds}{$j$-seeds}

\newcommand{\Uk}{U_k}

\title{Translation-Equivariant Matchings of  Coin-Flips on $\Z^d$}
\author{Terry Soo}

\address{Department of Mathematics, University of
 British Columbia, Vancouver BC, V6T 1Z2, CANADA}
\email{tsoo@math.ubc.ca}

\begin{document}
\begin{abstract}
Consider independent fair coin-flips at each site of the lattice $\Z^d$.  A translation-equivariant matching rule is a perfect matching of heads to tails that commutes with translations of $\Z^d$ and is  given by a deterministic function of the coin-flips.   Let $Z_{\Phi}$ be the distance from the origin to its partner, under the translation-equivariant matching rule $\Phi$.  Holroyd and Peres {\cite{Extra-Heads}} asked what is the optimal tail behaviour of $Z_{\Phi}$, for  translation-equivariant perfect matching rules.  We prove that for every $d \geq 2$, there exists a translation-equivariant perfect matching rule $\Phi$ such that $\E Z_{\Phi}^{\frac{2}{3} - \varepsilon} < \infty$ for every $\varepsilon > 0$.  
\end{abstract}
\maketitle
%
%\footnote{\hspace{-2em}Funded in part by the Natural Science and Engineering Research Council of Canada (NSERC PGS D) and a U.B.C. Graduate Fellowship  }
%\footnote{\hspace{-2em}{\bf Address:} {\tt

%\footnote{\hspace{-2em}{\bf Key words:} 
%\keywords{matching, Bernoulli random field}
%\footnote{\hspace{-2em}{\bf 2000 Mathematics Subject
%Classifications:} 
%\ams{60G55, 60G60}{60K35 }
%\end{abstract}
%
\section{Introduction}  
Consider the probability space $(\X, \F, \P)$, where $\X = \coin$, $\F$ is the standard product $\sigma$-algebra of $\coin$, and $\P$ is the product measure on $\F$ with parameter $p=\frac{1}{2}$.  We call elements of ${\Z^d}$ {\df{sites}}.
%and elements of $\{0,1\}^{\Z^d}$ {\df{configurations}}.  
%Let $\gamma$ be a configuration.  
%We say that a site $x \in \Z^d$ is {\df{occupied}} if $\gamma(x) = 1$ and {\df{unoccupied}} if $\gamma(x) = 0$.  
For $\gamma \in \X$, a  bijection $\phi:\Z^d \to \Z^d$ is a {\df{matching on $\gamma$}} if every site $x$ with $\gamma(x)=1$ is mapped to a site $y$ with $\gamma(y) =0$ and vice-versa and if the composition $\phi \circ \phi$ is the identity mapping on $\Z^d$.  For a site $z$, we define the translation $\theta^z$ on $\Z^d$ and $\X$ as follows: we set $\theta^z x := x+z$ for all $x \in \Z^d$ and for all $\gamma \in \X$ we set $\theta^z \gamma(x) := \gamma(x-z)$ for all $x \in \Z^d$.    A measurable mapping $\Phi: \coin \times \Z^d \to \Z^d$ is a {\df{matching rule}} if $\Phi(\gamma, \cdot)$ is a matching on $\gamma$ for $\P$-almost all $\gamma$.  We say that $\Phi$ is {\df{translation-equivariant}} if it commutes with translations; that is $\Phi(\theta^z \gamma, \cdot) = \theta^z \Phi(\gamma, \cdot)$ for $\P$-almost all $\gamma$.  

Let $\|\cdot \|$ be the $l^ {\infty}$ norm  on $\Z^d$.  We define $Z = Z_{\Phi}(\gamma):= \| \Phi(\gamma,\Oo) \|$ to be the distance from the origin $\Oo = (0, \dots, 0)\in \Z^d$ to its partner.  We will construct a translation-equivariant matching rule $\Phi$ and obtain upper bounds on $\PPp{Z > r}$.
\begin{theorem}
\label{result}
For all $d \geq 1,$ there exists a  translation-equivariant matching rule $\Phi$ such that for all $r > 0$, we have
\[  \result, \]
for some $c=c(d) < \infty $, where $\conditions$. 
\end{theorem}  
Prior to this result the best known decay appears to have been the following.
\begin{theorem}{\rm{  {\cite{meshalkin}},{\cite{Extra-Heads}}}}
\label{meshalkin}
For $d \geq 1,$ there exists a translation-equivariant matching rule $\Phi$ such that for all $r > 0$, we have  $\resultmm$,
for some $c=c(d) < \infty$. 
\end{theorem}
Theorem {\ref{meshalkin}} can be deduced from a simple construction due to Meshalkin.  Meshalkin matching was originally used to construct isomorphisms of Bernoulli schemes {\cite{meshalkin}}; it is the following construction.    In $d=1$, we define a translation-equivariant matching rule inductively, by first matching a zero to a one whenever a zero is immediately to the left of a one:
\[ \dots  {\bf\underline{01}}10 {\bf\underline{01}}1111000 {\bf\underline{01}}  \dots \]
In the next stage, we remove the matched pairs, and then follow the same procedure.  It is straightforward to check that bounding  $\PPp{Z > r}$ amounts to bounding $R = \inf \naiveset{m \geq 1: S_m = 0}$, where $S_m$ denotes a simple symmetric random walk. 

We may deduce the $d \geq 2$ case of Theorem {\ref{meshalkin}} from the following observation.   By applying a translation-equivariant matching rule $\Phi_{d-1}$ on $\Z^{d-1}$ to each $(d-1)$ dimensional {\em{plane}}, given by $\naiveset{z} \times \Z^{d-1}$ for each $z \in \Z$, we obtain a translation-equivariant matching rule $\Phi_{d}$ on $\Z^d = \Z \times \Z^{d-1}$, where $\PPp{Z_{\Phi_{d-1}} > r} = \PPp{Z_{\Phi_{d}} > r}$ for all $r > 0$. %It is easy to see that we do obtain a translation-equivariant matching rule and also $\PPp{Z_{\Phi_{d-1}} > r} = \PPp{Z_{\Phi_{d}} > r}$, for all $r > 0$. 

Theorem {\ref{result}} provides faster decay than that provided by Theorem {\ref{meshalkin}}, for all $d >1$.  After this paper was written, Tim\'ar {\cite{Timarb}} proved the following stronger result. 
\begin{theorem}{\rm  {\cite{Timarb}}  }
\label{cltconjecture}
For any $d \geq 1$  {there exists a translation-equivariant matching rule $\Phi$ such that for all $r > 0$,} we have $\resultalld$,
for some $c=c(d) < \infty$. 
\end{theorem}
\noindent Some of the methods of this article also appear in {\cite{Timarb}}.  New ideas are introduced in {\cite{Timarb}} and the methods of {\cite{Timarb}} are much more sophisticated.  

For $d=1,2,$ the bounds obtained in Theorem {\ref{cltconjecture}} are essentially best possible.  %The following theorem from {\cite{Extra-Heads}} provides a lower bound in $d=1,2$ on how fast $\PPp{Z_{\Phi} > r}$ can decay.
\begin{theorem}\rm{  \cite{Extra-Heads} }
\label{liggettbound}
If $d=1, 2$,  then for any translation-equivariant matching rule $\Phi$ we have $\E Z_{ \Phi} ^{\frac{d}{2}}  = \infty$. 
\end{theorem}
\noindent Hence for $d=1,2$ there does not exist a translation-equivariant matching rule $\Phi$ where  $\PPp{Z_{\Phi} > r} \leq c r^{-\rho}$, for some constants $\rho > \frac{d}{2}$ and $c = c(d) < \infty.$

Notice that the result of Theorem {\ref{liggettbound}} is only valid for $d=1,2$.  In fact for $d \geq 3$, Tim\'ar {\cite{Timarb}} has shown that one can find a translation-equivariant matching rule with even faster decay than that given by Theorem {\ref{cltconjecture}}.   
\begin{theorem}\rm{  \cite{Timarb} }
\label{expadam}
For any $d \geq 3$ and any $\e > 0,$  {there exists a translation-equivariant matching rule $\Phi$ such that for all $r > 0,$} we have $$\PPp{Z > r} \leq C\exp(-cr^{d-2-\e}),$$
for some constants $0<c,C  < \infty$. 
\end{theorem}

Variants of matching in continuum settings have also been studied; see {\cite{random}}, {\cite{Stable-PL}}, {\cite{Stable-Pltail}}, {\cite{Sodin}} and the references within.   

\paragraph{Outline of the Proof}  The proof of Theorem \ref{result} will proceed in two steps.  We will construct a translation-equivariant matching and then determine bounds for it.  To construct a translation-equivariant matching we will define, in a measurable translation-equivariant way,  a sequence $P_n$ of successively coarser partitions of $\Z^d$.  Following \cite{Hol-trees}, we call $P_n$ a {\df{clumping rule.}}  The members of $P_n$ are called {\df{clumps}} or {\df{$\boldsymbol{n}$-clumps}}, and we call the clumping rule {\df{locally finite}} if all the clumps are bounded.  A {\df{component}} of a clumping rule is a limit of some increasing (with respect to set inclusion) sequence of clumps.   A clumping rule is {\df{connected}} if it has only one component.   Adapting the construction in {\cite{Hol-trees}} we will construct a locally finite connected clumping rule.  From a  locally finite connected clumping rule it is easy to obtain a  translation-equivariant  matching rule $\Phi;$  this is because a translation-equivariant matching rule can be defined be first matching as many sites as possible within each $1$-clump and then iteratively matching as many unmatched sites as possible in each $n$-clump, for $n =2,3, \dots.$    We will obtain with the central limit theorem and a version of the mass transport principle, a preliminary result which implies that for $d \geq 3$  and for all $\e > 0$, we  have  $\PPp{Z_{\Phi} > r} \leq cr^{-\frac{3}{5} - \epsilon}$, for some constant $c = c(d, \e) < \infty.$  The preliminary result will not provide faster decay than that given by Theorem {\ref{meshalkin}} in the case $d=1,2$.  Upon a closer analysis of the geometry of the clumps, we will show that  clumps that are {\em{long and thin}} happen with small probability; this analysis is behind proof of  Theorem {\ref{result}}.

\paragraph{Outline of Paper}  The rest of the paper proceeds as follows.  In Section {\ref{sectiononclumps}} we discuss clumping rules and matchings from clumping rules.  In Section {\ref{seedsblobs}} we outline the construction of a clumping rule and collect some useful bounds.  In Section {\ref{masstransportsection}} we introduce a version of  the mass transport principle that will be useful in the proof of Theorem {\ref{result}}.  In Section {\ref{proof}} we prove Theorem {\ref{result}}.  We conclude the paper with some related open problems. 
{\nocite{firstLig}}
{\nocite{Hol-Lig}}
{\nocite{Thorissontrv}}
\section{Clumps}
\label{sectiononclumps}
Let $\Ppp_F(\Z^d)$ denote all finite subsets of $\Z^d$.  For $A \subset \Z^d$ define translations of $A$ via  $\theta^z A:= \naiveset{\theta^z x : x \in A}$.  Formally, a {\df{locally finite connected clumping rule}} is a measurable mapping $\CC:\coin \times \N \times \Z^d \to \Ppp_F(\Z^d)$ with the following properties.  For all $\gamma \in \coin,$ all $\  n \in \N$, and all  $x,y,z \in \Z^d$, we have
\begin{eqnarray*}
{\textrm{(i)}} && x \in  \CC(\gamma, n, x) \\
{\textrm{(ii)}} && \CC(\gamma, n, x) \cap \CC(\gamma, n, y) \not = \emptyset \Longrightarrow \CC(\gamma, n, x) = \CC(\gamma, n, y)  \\
{\textrm{(iii)}} &&  \CC(\gamma, n, x ) \subset \CC(\gamma, n+1, x) \\
{\textrm{(iv)}} &&  \CC(\theta^z \gamma, n, \theta^z x) = \theta^z \CC(\gamma, n, x) \\ 
{\textrm{(v)}} &&  \bigcup_n \CC(\gamma, n, \Oo) = \Z^d.
\end{eqnarray*}
Properties (i) and (ii) assure us that for each $n \in \N,$ the map $\CC(\cdot,n, \cdot)$ is a partition.  Property (iii) makes the partition successively coarser, (iv) is  translation-equivariance and (v) is connectedness.  
\begin{proposition}{There exists a locally finite connected clumping rule almost surely.}
\label{existsc}
\end{proposition}
\noindent The proof of Proposition {\ref{existsc}} will be given in the next section.
\paragraph{Matchings from Clumpings.}  From a locally finite connected clumping rule, we can   construct a translation-equivariant matching rule in a countable number of stages.  In the first stage, within each of the $1$-clumps we match every possible site.   Given that the $(n-1)$ stage is completed, within each of the $n$-clumps we match every site we can, ignoring the sites that were previously matched.  In order to ensure that the resulting matching is translation-equivariant, use for example a lexicographic ordering on $\Z^d$, to determine the maximal partial matching on the clumps.   Ergodicity, connectedness, and the fact that $p=\frac{1}{2}$ gives us that every site will be matched at some stage.  Note that for our purposes we do not {\em{need}} to make this argument as we obtain upper bounds on $\PPp{Z > r}$ which easily imply that $\PPp{Z > r} \to 0$, as $r \to \infty$, (see Theorem {\ref{result}} or Proposition {\ref{firstresults}}).

In the next section, we construct an explicit locally finite connected clumping rule $\CC$.   %In fact, any translation-equivariant matching rule derived from this clumping using a construction such as the one described above will satisfy the conditions of Theorem {\ref{result}}. 

\section{Seeds, Cutters and Blobs}
\label{seedsblobs}
Our construction of the clumping rule $\CC$ is adapted from  {\cite{Hol-trees}}.  In {\cite{Hol-trees}} and {\cite{Timar-trees}} clumpings are used  to obtain factor graphs of point processes.  See also {\cite{{MR2044812}} for background.
\subsection{Basic Set-Up}
Let $\|\cdot \|$ denote the $l^ {\infty}$ norm on $\Z^d$.     Let  $S(x,r):=\{y \in \R^d: \|x-y\| \leq r\}$.  Thus $S(x,r)$ is the cube of side length $2r$ centered at $x$.  We also write $S(\Oo,r) = S(r)$.  We let $\naiveset{e_m}_{m=1} ^d$  be the standard unit basis vectors in $\R^d$.

For each $k \in \N$, we say that a site $x\in \Z^d$ is a {\df{$\boldsymbol{k}$-seed}} if $\gamma(x) =1$, and 
 $\gamma(y) = 0 \textt{for all}  y \in \{x + n e_1 : 1 \leq n \leq k-1 \}$.
Whenever $x$ is a \kseed \  we call the set $\{x + n e_1 : 0 \leq n \leq k-1 \}$ its {\df{shell}}.
For example, a $4$-shell has the form:
\[ {\bf{1}}000. \]
Note that the probability of a \kseed \  occurring at a particular point is exactly $2^{-k}$.  
%\begin{equation}
%\label{akk}
%{   
%a_k:= \obracket{\frac{1}{2}}^{{k}}.  
%}
%\end{equation}
%
Two seeds are said to be {\df{overlapping}} if their shells intersect.  %Although two seeds can overlap, for a fixed $k$, no two distinct \kseeds \ can  overlap.   
Note that two seeds $x$ and $y$ overlap if and only if $x = y$.  
This property will be useful later (see Section {\ref{secondattemptsection}}).
We define
\begin{equation}
\label{rk}
{
r_k =  {  {  \big( 2^k{k^2}  \big)   }^{\frac{1}{d}} } +\frac{1}{2}. 
}
\end{equation}
The reason for the choice of $r_k$ will be evident shortly.  %The addition of the term $\frac{1}{2}$ is to ensure that $r_k$ is never an integer.  
Define the vector 
\begin{equation}
\label{sssk}
{
s_k: = \floor{100{r_k}}e_1.
}  
\end{equation}
A {\df{$\boldsymbol{k}$-cutter}} is a subset of $\R^d$ of the form $\naiveset{ y \in \R^d: \|y-x\| = r_k }$,
where $x - s_k$ is a \kseed.  We introduce a shift $s_k$ for technical reasons which will surface latter.  
We define $W_k \subset \R^d$ to be the union of all the \kcutters.  Note that we have chosen $r_k$ so that $r_k \not \in \N$.  Thus we have that $W_k \cap \Z^d = \emptyset$ for all $k \in \N$.      A {\df{$\boldsymbol{k}$-blob}} is a connected component of $\R^d \setminus  \cup_{j > k} W_j$. 
Hence we have that the sequence of \kblobs \ define a successively coarser partition of $\R^d$ (ignoring the elements of $\cup_k W_k$.)   The \kblobs \ induce a clumping rule ${\CC}$ when we intersect the \kblobs \ with $\Z^d$.  Note the technical distinction between blobs and clumps.  

It is obvious that the induced clumping rule ${\CC}$ is translation-equivariant; it remains to show that it is locally finite and connected.  It suffices to show that all the blobs are bounded and that for every $x \in \R^d$, there is a \kblob \ that contains both $x$ and the origin.

%In order to show that the clumping rule $\CC$ is locally finite it suffices to show that the blobs are bounded.  In order to show that the clumping rule $\CC$ is connected, it suffices to show  that for every $x \in \R^d$ there is a \kblob \ that contains both $x$ and the origin.  
%
%
%
\subsection{Estimates}
\label{estimates}
In this section we obtain some estimates that will show that the clumping rule $\CC$  as defined in the previous section is indeed locally finite and connected.   The following events will be important in our analysis.  Let 
\begin{equation}
\label{disek}
E_k(x):= \naiveset{ x \textt{is enclosed by some} k\text{-cutter} };
\end{equation} 
that is $E_k(x)$ occurs if and only if  for some site $x_0$, $x_0 - s_k$ is a $k$-seed and \\ $x \in \naiveset{ y \in \Z^d : \| y - x_0 \| \leq r_k }$. Also  let $E_k = E_k(\Oo)$.  Let $$\Uk(s):= \naiveset{S(s) \textte{intersects some $k$-cutter}};$$ that is $\Uk(s)$ occurs if and only if for some site $x_0$, $x_0 - s_k$ is a $k$-seed and \\  $\naiveset{ y \in \R^d : \| y - x_0 \| = r_k } \cap S(s) \not = \emptyset$.  Also let 
\begin{equation}
\label{disck}
C_k(s): = \cup_{j \geq k} \Uk(s).  
\end{equation}
From an analysis of these events, we will deduce that the clumping rule $\CC$ is both locally finite and connected.    Moreover, we will see that the tail behaviour of $Z_{\Phi}$ (where $\Phi$ is a translation-equivariant matching rule obtained from the clumping rule $\CC$) also depends on these events.
%
%

%When the occasion arises we will denote positive finite constants that possibly depend on $d$, the dimension, using $c_1, c_2, \dots$.
%
%
%
\begin{lemma}[Enclosure bounds]
\label{ENC}
For all $k > \cconstantk$, for some $\cconstantk= \cconstantk(d) < \infty$, we have  $\PPp{E_k ^c} \leq  e^{-k  }$. 
\end{lemma}
\begin{proof}
Note that, 
\begin{equation}
\label{VKK}
{
E_k =  \naiveset{ S(-s_{k}, r_{k} - 1) \textte{contains some \kseed}}.  
}
\end{equation}
Let $p_k$ the maximum possible number of \kseeds \ inside $S(-s_{k}, r_{k} - 1)$.    Recall that no two (distinct) $k$-seeds overlap and the probability that a $k$-seed occurs at a particular point is $2^{-k}$.  Hence, $\PPp{E_k} \geq 1 - (1 - 2^{-k})^{p_{k}} \geq 1 - e^{-2^{-k} p_k}.$
By our choice of $r_k$ in  ({\ref{rk}}) and since $(\frac{r_k^d}{k}) \leq  p_k \leq  \lceil \frac{2r_k}{k} \rceil (2r_k)^{d-1} $ for all $k \geq \cconstantk$, for some $\cconstantk =\cconstantk(d) < \infty$, we have that $\PPp{E_k} \geq 1 - e^{-k}$.  
\end{proof}
%\noindent {\em{Remark}}: Note that by translation-equivariance the same result applies to $E_k(x)$ for any site $x.$
%
%
\begin{corollary}{All ${k}$-blobs  are bounded almost surely.}
\label{bbound}
\end{corollary}

\begin{proof}
It suffices to show that all \kblobs \  that contain $\Oo$ are bounded.  By Lemma {\ref{ENC}}, we have that $\P(E_k) \to 1$ as $k \to \infty$, so that $E_k$ occurs for infinitely many $k$ almost surely.  Hence all blobs which contain $\Oo$ are bounded.      
\end{proof}
\begin{lemma}[Cutter bounds]
\label{CUT}
For all $ k \geq 1$ and all $s > 0$, we have  $\PPp{C_k(s)}  \leq \ctimesfirstterm s  \big( \frac{k^{ \old }}{r_{k}} \big)$, 
for some $\ctimesfirstterm = \ctimesfirstterm(d) > 0$.
\end{lemma}
\begin{proof}
Observe that
\begin{equation}
\label{AKK}
{
\Uk(s) =  \naiveset{S(-s_k, r_k + s) \setminus S(-s_k, r_k - s) \textte{contains some \kseed}}.  
}
\end{equation}
Clearly,  $ \PPp{\Uk(s)} \leq N_k(s)2^{-k}$, where  $N_k(s)$ is the number of lattice points in  $S(-s_k,r_k + s) \setminus S(-s_k,r_k - s)  $.  We have that $ N_k(s) = |S(r_k +s)| - |S(r_k - s)| \leq \clatticep r_k^{d-1} s $, for some $\clatticep = \clatticep(d) > 0$.  So we obtain that, 
$\PPp{\Uk(s)} \leq \clatticep s r_k ^{d-1} 2^{-k}$.
Thus recalling our choice of $r_k$ in ({\ref{rk}}), we have that
\begin{equation}
\label{evident}
{
\PPp{C_k(s)} \leq \sum_{ j \geq k} \P(U_j(s)) \leq \clatticep s \sum_{j \geq k} 2^{-k} r_k^{d-1} \leq \ctimesfirstterm s  \bigg( \frac{k^{ \old  }}{r_{k}} \bigg).
}
\end{equation}
\end{proof}
\begin{corollary}{The clumping rule $\CC$ is connected almost surely.}
\label{cconn}
\end{corollary}

\begin{proof}
Let $s > 0$.  By the Borel-Cantelli lemma and (\ref{evident}) we have that  $\Uk(s)$ occurs infinitely often with probability zero.  Thus any point within distance $s$ of $\Oo$ will eventually share a blob with it.   
%
%Consider an arbitrary $s > 0$.  We show that $\Uk(s)$ occurs infinitely often with probability %zero.  Hence any point within distance $s$ of $\Oo$ will eventually share a blob with it.   By% the Borel-Cantelli lemma it suffices to show that $\sum_{k=1} ^ { \infty } \P(\Uk(s)) < %\infty,$  which follows from (\ref{evident}).
%
%
\end{proof}
\begin{proof}[Proof of Proposition {\ref{existsc}}]
Apply Corollaries {\ref{bbound}} and {\ref{cconn}}.
\end{proof}

Now we obtain  a translation-equivariant matching rule $\Phi$ from our locally finite connected clumping rule $\CC$, via the procedure outlined in Section {\ref{sectiononclumps}}.   We will use Lemmas {\ref{ENC}}, {\ref{CUT}, and the central limit theorem to obtain bounds on $Z_{\Phi}$.  %(Recall $Z_{\Phi}$ is the distance from the origin to its partner under the translation-equivariant matching rule $\Phi$.)
\section{Mass Transport}
\label{masstransportsection}
We will require a version of the mass transport principle in order to facilitate calculations.  See {\cite{two}} and {\cite{seven}} for background.  Our main application of the mass transport principle  will be to prove a modified version of Lemma {\ref{MASSeg}} below, which states that each site has an equal chance of not being matched within its $k$-clump.  Similar ideas also appear in {\cite{random}}.   

Let $\CC$ be the clumping rule defined in Section {\ref{seedsblobs}} and let $\Phi$ be the corresponding translation-equivariant matching rule obtained from $\CC.$  We say that a site is {\df{$\boldsymbol{k}$-bad}} if it is not matched in its $k$-clump.  Let $\clump_k(x)$ be the $k$-clump containing the site $x$ and let $\clump_k(\Oo) = \clump_k$ be the $k$-clump containing the origin.  Let $\cardsupp {\clump_k }$ be the cardinality of $\clump_k$.  Consider the sum
$\zeta:=\sum_{x \in \clump_k} (2 \gamma(x) -1)$, so that $|\zeta|$ is the number of $k$-bad  sites in $\clump_k$.  %(Recall that we match every site that we can at each stage.)   
%
%\begin{equation}
% \zeta(\clump_k(y)):= \sum_{x \in \clump_k(y)} (2 \gamma(x) -1),  
%\end{equation}
%
%so that $|\zeta(\clump_k)| = \textt{number of} k \text{-bad sites in} \ \clump_k$.  (Recall that we match every site that we can at each stage.)   
%
%
%
\begin{lemma}
\label{MASSeg}
For all $k \geq 1,$ the probability that the origin is k-bad is exactly
  \[  \E \bigg( \frac{1}{\cardsupp{\clump_k}}\bigg| \sum_{x \in \clump_k} (2 \gamma(x) -1) \bigg |  \bigg). \]
\end{lemma}

We define a {\df{mass transport}} to be a non-negative measurable function $T: \Z^d \times \Z^d \times \coin \to \R$ which is translation-equivariant; that is for all $x,y \in \Z^d$ and for all $\gamma \in \coin$ and for all translations $\theta$ of $\Z^d$, we have $T(\theta x, \theta y, \theta \gamma) = T(x,y,\gamma).$
For $A, B \subset \Z^d$, we let $T(A, B, \gamma):=  \sum_{x \in A, y \in B} T(x,y, \gamma).$
We think of $T(A,B,\gamma)$ as the mass transfered from $A$ to $B$ under $\gamma \in \X$.  %Recall that $\theta^y$ is the  translation on $\Z^d$ where $\theta^y x = x+ y$.  
We will use the following version of the mass transport principle.   %An application of the Fubini theorem and translation-equivariance of $T$ and translation-invariance of $\P$ gives us the following.
\begin{lemma}[Mass transport principle]
\label{MASS}
For any mass transport, $T:\Z^d \times \Z^d \times \coin \to \R$, we have $\E T(\Oo, \Z^d, \cdot) = \E T(\Z^d, \Oo, \cdot)$.
\end{lemma}
\begin{proof}
\begin{eqnarray*}
 \E T(\Oo, \Z^d, \cdot) &=& %\int \sum_{y \in \Z^d} T(\Oo, y, \gamma )d\P(\gamma) 
\sum_{y \in \Z^d }  \int T(\Oo, y, \gamma)d\P(\gamma) \\ 
&=& \sum_{y \in \Z ^d} \int T(-y, \Oo, \theta^{-y} \gamma)d\P(\gamma) \\
&=& \sum_{y \in \Z ^d} \int T(-y, \Oo, \gamma)d\P(\gamma) \\
&=& \E T(\Z^d, \Oo, \cdot).
\end{eqnarray*}
The first and last equalities follows from the Fubini theorem.  The second equality follows from the translation-equivariance of $T$ and the third equality follows from the translation-invariance of $\P$.  
\end{proof}
To illustrate the versatility of the mass transport principle (Lemma {\ref{MASS}})  we prove the following (unsurprising) fact.
\begin{proposition}
\label{nomatchinghalf}
Let $\F$ be the standard product $\sigma$-algebra of $\coin$ and let $\P_p$ be the product measure on $\F$ with parameter $p$.   If $p \not = \frac{1}{2}$, then there does not exists a translation-equivariant matching rule. 
\end{proposition}
\begin{proof}%[Proof of Proposition {\ref{nomatchinghalf}}]
Let $\Phi$ be a translation-equivariant matching rule.  Consider the mass transport $M$ defined as follows.  Let $x$ be a site and $\gamma \in \X$.  If $\gamma(x)=1$, then $M(x, x, \gamma)=1$; that is $x$ sends one unit of mass to itself.  Otherwise $M(x, y, \gamma) =1$, where $y$ is a site with $\Phi(x,\gamma)=y$ and $\gamma(y)=1$; that is $x$ sends out a unit of mass to the site $y$ that it is matched to under $\Phi(\cdot, \gamma)$.  Since $\Phi$ is translation-equivariant this defines a mass transport $\P_p$-a.s.  Let $\E_p$ be the expected value operator with respect to the measure $\P_p$.   Now since every site sends out exactly one unit of mass we have that $\E_p M(\Oo, \Z^d, \cdot) = 1$.
Also by considering the cases $\gamma(\Oo) =1$ or $\gamma(\Oo) =0$ , we also have that $\E_p M(\Z^d, \Oo, \cdot) = 2p $.
Hence we  have by the mass transport principle that $p= \frac{1}{2}$.
\end{proof}

\begin{proof}[Proof of Lemma {\ref{MASSeg}}]
For each $k \geq 1,$ we define a  mass transport $T_k$ by saying, if a site $x$ is $k$-bad, then it sends out one unit of mass uniformly to every site in its $k$-clump $\clump_k(x)$, while  $x$ sends out no mass, if it is not bad.  To be precise,
\begin{equation}
\label{tkkk}
T_k(x,y, \gamma): = \left(\frac{1}{ \cardsupp{\clump_k(x)}} \right) \one{ x \textte{is \kbad} }(\gamma) \one{y \in L_k(x) } (\gamma).
\end{equation}
It is easy to see that  
%
%\begin{eqnarray*}
$$ \E T_k(\Z^d, \Oo, \cdot) = \E \bigg( \frac{1}{\cardsupp{\clump_k} }\sum_{x \in \clump_k} \one{ x \textte{is \kbad} } \bigg).$$  
%&=& \E \big( |\zeta(\clump_k) | / \cardsupp{\clump_k} \big).
%\end{eqnarray*}
On the other hand, we have  that $\E T_k(\Oo, \Z^d, \cdot) = \PP{ \Oo \textte{is \kbad}}$.
Thus, an application of  mass transport principle completes the proof.
\end{proof}

Now we are in a position get bounds on $\PPp{Z > r}$.  We will see that mass transport principle with information about the size of $L_k$ and its diameter gives us an estimate with an application of the central limit theorem.      
\section{Proof of Theorem {\ref{result}} }
\label{proof}
\subsection{First Estimates}
\label{firstattemptsection}
Let $\Phi$ be the translation-equivariant matching rule we obtain from the clumping rule ${\mathcal C}$ defined in Section {\ref{seedsblobs}}.   Recall that $Z= Z_{\Phi}$ is the distance from the origin to its partner under $\Phi$.    We will obtain bounds on $\PPp{Z > r}$ by choosing a sequence of events $D_k$ and a $K=K(r)$ so that $\naiveset{Z > r} \cap D_K \subset \naiveset{\Oo \ \text{is \Kbad} }$.  The events $D_k$ will be chosen in a way so that we can obtain upper bounds on  $\PP{\Oo \ \text{is \Kbad} }$ and $\PPp{D_K ^c}.$

Let $\a \in (0,1)$.  In fact we will end up choosing  $\a = \a(d) = \choiceofalpha$.  Recall that the events $E_k$  and $C_k(s)$ were defined earlier in Section \ref{estimates}; see  (\ref{disek}) and (\ref{disck}).  Let $B_k$ be the \kblob \ containing the origin.  The following relations describe the geometry of $B_k$, when $E_k$ or $C_k(r^ \a)$ occur.  We have
\begin{equation}
\label{TERME}
{
 E_k \subset \naiveset{ \text{there exists} \ x \  \text{so that} \  B_k \subset S(x,r_{k}) \subset S(2r_{k}) }
}
\end{equation}
%
%\begin{equation}
%\label{TERME}
%{
%    E_k  \implies \text{there exists} \ x \  \text{so that} \  B_k \subset S(x,%r_{k}) \subset S(2r_{k})
%}
%\end{equation}
%
and
\begin{equation}
\label{TERMC}
{
 C_k (r^ \a)^c \subset \naiveset{S(r^\a) \subset B_k}. 
}
\end{equation}
%\begin{equation}
%\label{TERMC}
%{
% C_k (r^ \a)^c  \implies S(r^\a) \subset B_k.
%}
%\end{equation}
%    
We consider the following decomposition: 
\begin{equation}
\label{TERMS}
{ 
\naiveset{ Z > r } \hspace{0.2 cm} = \hspace{0.2 cm} \Big(  \obracket{E_k \cap C_k(r^ \a)^c }    \cap \naiveset{ Z > r  }   \Big) \hspace{0.1 cm} \cup \hspace{0.1 cm}   \Big( \obracket{E_k^c \cup C_k(r ^ \a)}  \cap \naiveset{ Z > r } \Big) .  
} 
\end{equation}
See Figure {\ref{ekckevent}}  for a realization of the event $E_k \cap C_k(r ^\a)^c$.
\begin{figure}[htbp]
\begin{center}
\includegraphics[width=3.5cm]{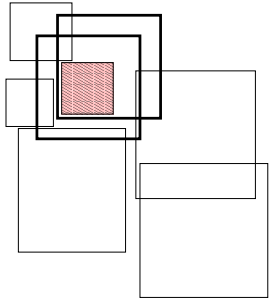}
\end{center}
\caption{An illustration of the event $E_k \cap C_k( r^\a)^c$.  The thick cutters represent $k$-cutters enclosing the origin.  This corresponds to the event $E_k$.   The shaded region represents the no cutter zone of radius $r^\a$ about the origin.  This corresponds to the event $C_k(r^\a)^c$.  }
\label{ekckevent}
\end{figure}

The role of the parameter $\a$ can be explained heuristically as follows.  If the parameter  $\a$ is small, then $C_k(r^\a)^c$ occurs with high probability, but then $B_k$ could possibly be very small.  If $\a$ is close to $1$, then $B_k$ would almost contain a cube of length $2r$, but then $C_k(r^\a) ^c$ occurs with low probability.  We will choose $\a$ to optimize over these alternatives.

Let $K=K(r)$ be defined to be the unique integer $K$ such that

\begin{equation}
\label{choiceofk}
{   
2r_K < r_{K+1} < r \leq r_{K+2}  
}
\end{equation}
Note that for some $\cklog = \cklog(d) > 0$, we have that for all $k \geq \cklog$,  $r_k = (2^k k^2)^{\frac{1}{d}} + \frac{1}{2} \leq (\frac{e^k}{2^{-k}})^{\frac{1}{d}} $.  Hence applying ({\ref{rk}}) with ({\ref{choiceofk}}) for all $r$ sufficiently large we have that
\begin{equation}
\label{kdbound}
{
  \Big(\frac{\log2}{1 + \log2} \Big) (K+1)   \leq \frac{d \log r}{ 1 + \log 2} \leq K+2.  
}
\end{equation}
\begin{proposition}[Decay of the first term in ({\ref{TERMS}})  via the central limit theorem]
\label{CLTdecay}
For all $r > 0$ and for  the unique integer $K=K(r)$ such that $r_{K+1} < r \leq r_{K+2}$, we have
\[ \P \Big(  \obracket{E_K  \cap C_K(r^ \a)^c } \cap \naiveset{ Z > r  } \Big) \leq  \frac{\cclt}{ (r^{\a})^{  \fracc{d}{2}  } }, \]

for some $\cclt = \cclt(d) > 0.$
\end{proposition}
\noindent Remark: Note that the decay in Proposition {\ref{CLTdecay}} is the decay that appears in Theorem {\ref{cltconjecture}}, if we set $\a =1$.

Before we begin the proof of Proposition \ref{CLTdecay}, we collect some easy, but important observations.   By ({\ref{TERME}})  and ({\ref{choiceofk}}), we have
\[  \obracket{E_K  \cap C_K (r^ \a) ^c  } \cap \naiveset{Z > r} \ \subset \  \naiveset{\Oo \ \text{is \Kbad} }.\]
So from ({\ref{TERMC}}) we have 
\begin{equation}
\label{condind}
\big(  \obracket{E_K  \cap C_K(r^ \a)^c } \cap  \naiveset{ Z > r  } \big) \ \subset \    \big( \naiveset{\Oo \ \text{is \Kbad}}   \cap  {E_K  \cap \naiveset{ \cardsupp{\clump_K} \geq r^{\a d}} }  \big).
\end{equation}
Recall that $\clump_k$ is the $k$-clump containing $\Oo$.  To analyze the right hand side of  ({\ref{condind}}) we will use the following version of Lemma {\ref{MASSeg}}. 
\begin{lemma}
\label{lemmaaftermass}
For all $ k \geq 1$, we have
\begin{eqnarray}
\label{aftermass}
&\hspace{0.1 cm}& \P \Big( \naiveset{\Oo \ \text{is $k$-{\rm{bad}}}}  \cap  E_k \cap \naiveset{\cardsupp{\clump_k\geq r^{\a d}} } \Big) =  \nonumber \\
&\hspace{0.1 cm}& \E \Bigg( \frac{1}{\cardsupp{\clump_k}}\bigg| \sum_{x \in \clump_k} (2 \gamma(x) -1) \bigg| \ ; \ E_k \cap \naiveset{\cardsupp{L_k} \geq r^{\a d}}  \Bigg). 
\end{eqnarray}
\end{lemma}
\begin{proof}
We use the mass transport principle.  Recall $T_k$ as defined in the proof of Lemma {{\ref{MASSeg}}}: 
\[ T_k(x,y, \gamma): = \left(\frac{1}{ \cardsupp{\clump_k(x)}} \right) \one{ x \textte{is \kbad} }(\gamma) \one{y \in L_k(x) } (\gamma). \]
Define another mass transport: $$\hat{T}_k(x,y,\gamma) := T_k(x,y,\gamma){\one{ E_k(x) \cap  \naiveset{\cardsupp{\clump_k(x)} \geq r^{\a d} }}(\gamma)}.$$
Note that on the event $\naiveset{y \in L_k(x)}$ we have that the event $E_k(x)$ occurs if and only if the event $E_k(y)$ occurs and $\cardsupp{L_k(x)}= \cardsupp{L_k(y)}.$
Hence we obtain that 
\begin{eqnarray*}
 \E \hat{T}_k(\Z^d, \Oo, \cdot) &=&  \E \sum_{y \in \Z^d}T_k(y, \Oo, \cdot){\one{ E_k(y) \cap  \naiveset{\cardsupp{\clump_k(y)} \geq r^{\a d} }} } \\
&=& \E \sum_{y \in \Z^d}T_k(y, \Oo,\cdot){{\one{ E_k \cap  \naiveset{\cardsupp{\clump_k} \geq r^{\a d} }} }} \\
%&=& \E \Bigg( \frac{| \zeta(L_k)| }{\cardsupp{L_k}}; E_k \cap \naiveset{\cardsupp{L_k} \geq r^{\a d}}           \Bigg).
&=& \E \Bigg( \frac{1}{\cardsupp{\clump_k}}\bigg| \sum_{x \in \clump_k} (2 \gamma(x) -1) \bigg| \ ; \ E_k \cap \naiveset{\cardsupp{L_k} \geq r^{\a d}}  \Bigg). 
\end{eqnarray*} 
On the other hand, we have that
\begin{eqnarray*}
\E \hat{T}_k(\Oo, \Z^d, \cdot) &=& \P \Big( \naiveset{\Oo \ \text{is \kbad}}  \cap  {E_k  \cap \naiveset{ \cardsupp{\clump_k} \geq r^{\a d}} } \Big). 
\end{eqnarray*}
Thus an application of the mass transport principle (Lemma {\ref{MASS}}) completes the proof.    
\end{proof}
Next we will use the central limit theorem to estimate the right hand side of  ({\ref{aftermass}}), but first we need to verify that we have the necessary independence.  For $k \geq 1,$ consider the events:  $$H_k(x_1, x_2, \dots x_n):= \big\{ \naiveset{x_1, \dots, x_n}  =   L_k \cap S(2r_k)  \big \},$$ where $x_i \in \Z^d$ and  $\|x_i \| \leq 2r_k$.    Let $\G_k := \sigma(\gamma(x) : x \in \bgset)$.  The following lemma is behind why the (large) shift  $s_k = \floor{100{r_k}}e_1$, along the axis $e_1$, appears in the definition of the $k$-cutters.
\begin{lemma}
\label{independent}
For all $k \geq 1$ and for all $\|x_i\| \leq 2r_k,$ the $\sigma$-field  $\G_k$ is independent of $\sigma(H_k(x_1, \dots x_n), E_k)$.  
\end{lemma}
\begin{proof}
Consider a site $y$ with $\|y\| < \frac{s_k}{3}$.  The event $\naiveset{ y \in \clump_k }$ is determined by whether there exists $j$-seeds, with $j \geq k$, to give rise to $j$-cutters that can separate $y$ from $\Oo$.  However such \jseeds \  (and their shells) are at least at distance $\frac{s_k}{2}$ from the origin.  So we have that  $\naiveset{\gamma(x) : \|x\| < \frac{s_k}{3} }$ is independent of $H_k(x_1, x_2, \dots, x_n),$ for all $x_i \in \Z^d$ such that $\|x_i\| \leq 2r_k.$    Also recall that $E_k$  from ({\ref{VKK}}) is determined by $\gamma(x)$, where $\|x\| \geq \frac{s_k}{2}$.    \end{proof}
Now the proof of Proposition {\ref{CLTdecay}} amounts to a simple calculation, whose result we record in the next lemma.
\begin{lemma}
\label{againuse}
For all $k \geq 1$, we have 
\begin{equation*}
\P\Big( \naiveset{\Oo \ \text{is  $k$-{\rm{bad}}}}  \cap  {E_k  \cap \naiveset{ \cardsupp{\clump_k} \geq r^{\a d}} } \Big)  \leq  \fracc{\cclt}{ (r^{\a})^{  \fracc{d}{2}  }}, 
\end{equation*}
for some $\cclt > 0.$  
\end{lemma}
\begin{proof}[Proof of Proposition {\ref{CLTdecay}}]
From (\ref{condind}) and Lemma {\ref{againuse}},  Proposition {\ref{CLTdecay}} follows immediately.
\end{proof}  
\begin{proof}[Proof of {Lemma \ref{againuse}}]
Let $k \geq 1,$ recall that by ({\ref{TERME}}) we know that on the event $E_k$, we have $L_k \subset \bgset$.   Fix $x_1, \dots, x_n \in S(2r_{k})$ and let $H_k = H_k(x_1, \dots, x_n).$   We will now compute
\[ A:=  \E \Bigg( \frac{1}{\cardsupp{\clump_k}}\bigg| \sum_{x \in \clump_k} (2 \gamma(x) -1) \bigg| \ ; \ E_k \cap  H_k(x_1, \dots, x_n)           \Bigg), \]
by conditioning on $\G_k$.  Let  $S_n = \sum_{i=1}^n (2\gamma({x_i}) -1).$ 
Consider the following calculation:
\begin{eqnarray}
A  %&=& \E \Bigg(  \cardsupp{L_k}^{-1} \abs{ \sum_{x \in \Z^d} (2\gamma(x) -1) \one{ x \in L_k} }  \one{ E_k \cap H_k} \Bigg) \\
 & =&  \E \Bigg(  n^{-1} \abs{ \sum_{i=1}^n (2\gamma({x_i}) -1) } \one{ E_k \cap H_k} \Bigg) \\
\label{third}
& = & \E  \Big(  \E \Big( n^{-1} |S_n| \one{ E_k \cap H_n } \ \boldsymbol{\Big|} \ \G_k  \Big) \Big) \ \  \\
\label{fourth}
&=&  \E \Bigg( \frac{1}{n} |S_n|  \E \Big(  \one{E_k \cap H_k} \ \Big|  \  \G_k   \Big) \Bigg) \ \  \\
\label{fifth}
&=& \E \Big( \frac{1}{n} |S_n| \Big)  \E (  \one{ E_k \cap H_k}  ) \ \   \\
%&=&  \frac{\E|S_n|}{n} \PPp{ E_k \cap H_n}  \\
\label{final}
&\leq& \frac{\cclt}{\sqrt{n}} \PPp{ E_k \cap H_k}, \ \  
\end{eqnarray}
for some $\cclt > 0.$   Equality ({\ref{third}}) is obtained by conditioning on the $\sigma$-field $\G$.  Equality ({\ref{fourth}}) comes from the fact that the $\gamma(x_i)$ are all $\G$ measurable.  By Lemma {\ref{independent}}, we have that $\G_k$ and $\sigma(H_k, E_k)$ are independent, thus equality ({\ref{fifth}}) follows.  Inequality (\ref{final}) is obtained by applying the central limit theorem.  

By summing over all possible $H_k(x_1, \dots, x_n)$, we obtain that
\[  \E \Bigg( \frac{1}{\cardsupp{\clump_k}}\bigg| \sum_{x \in \clump_k} (2 \gamma(x) -1) \bigg|  \ ; \  E_k \cap \naiveset{\cardsupp{L_k}=n}           \Bigg) \leq \frac{\cclt}{\sqrt{n}} \PP{ \cardsupp{L_k} = n} . \]
Furthermore,  by summing over all $n \geq r^{\a d}$ we see that
\begin{eqnarray*}
%\label{againuse}
\E \Bigg( \frac{1}{\cardsupp{\clump_k}}\bigg| \sum_{x \in \clump_k} (2 \gamma(x) -1) \bigg|  \ ; \ E_k \cap \naiveset{\cardsupp{L_k} \geq r^{\a d}}           \Bigg) & \leq&  \fracc{\cclt}{ (r^{\a})^{  \fracc{d}{2}  } }. 
\end{eqnarray*}
Thus an application of Lemma {\ref{lemmaaftermass}} completes the proof.
\end{proof}
Now we turn our attention to the second term in ({\ref{TERMS}}): $\obracket{E_k^c \cup C_k(r ^ \a)} \cap \naiveset{ Z > r}$.
We will bound this term in two different ways.  As a first step, let us just throw away the term $\{Z > r\},$ since this will allow us to obtain a novel result for the case $d \geq 3$ without much more additional effort.   
\begin{lemma}[Decay of the second term in ({\ref{TERMS}}):  First bound]
\label{throwawaydecay}
For all $ r > 0$ and for the unique integer $K=K(r)$ such that $r_{K+1} < r \leq r_{K+2}$, we have
\[  \PPp{{E_K^c \cup C_K(r ^ \a)}} \leq    \cmethod1 r^{\a-1} (\log r)^{\old},  \]
for some $\cmethod1 = \cmethod1(d) > 0$. 
\end{lemma}
\begin{proof}
Recall that from Lemma {\ref{ENC}} and Lemma {\ref{CUT}}, we already have bounds for the events appearing in this term.  From ({\ref{kdbound}}) we see that %the $ \PPp{E_K^c}$ is not going to give us any problem; 
\begin{equation}
\label{encdecay}
{  \PPp{E_K^c} \leq \cestager r^{-\frac{d}{1 + \log 2}}, }
\end{equation}
for some $\cestager = \cestager(d) > 0$.
Note that $\frac{d}{1 +\log 2} > \frac{d}{2}$.  On the other hand, applying ({\ref{kdbound}}) to Lemma {\ref{CUT}} we obtain that
\begin{equation}
\label{firsttermrr}
{
 \PPp{ C_K(r^\a) } \leq \ccutterconst r^{\a -1} (\log r)^{\old},
}
\end{equation}
for some $\ccutterconst = \ccutterconst(d) > 0$. 
\end{proof}
\begin{proposition}[Easy preliminary result]
\label{firstresults}
For all $d \geq 1,$ there exists a translation-equivariant matching rule $\Phi$ such that $Z = Z_{\Phi}(\gamma) = \|\Phi(\gamma, \Oo)\|$ has the tail behaviour: $\P(Z > r) \leq \constantmethod1 r^{-\beta'} (\log r)^{\old}$, 
where $\constantmethod1= \constantmethod1(d) > 0$ and $\beta'= \beta'(d) = \frac{1}{1 + \frac{2}{d}}$.
\end{proposition}
\begin{proof}
We can see from ({\ref{TERMS}}) and Proposition {\ref{CLTdecay}} and Lemma {\ref{throwawaydecay}} that 
\[ \P(Z > r) \leq \cclt r^{\frac{-\a d}{2}} + \cmethod1 r^{\a -1} (\log r)^{\old}.\]  
Hence we are led to minimize the quantity:  $\max( \frac{-\a d}{2},    \a - 1)$.
So we choose (for the purposes of this proposition)  $\a = \a(d)= \frac{1}{1 + \frac{d}{2}}$.
\end{proof}
Proposition {\ref{firstresults}}, gives for $d=2$, decay of order $\frac{(\log r)^2}{r^{\fracc{1}{2}}}$.  For $d=2$, Theorem {\ref{meshalkin}} still provides a better result, but for $d \geq 3,$ Proposition 17 provides faster decay than Theorem {\ref{meshalkin}}.  
\subsection{Long and Thin Blobs}
\label{secondattemptsection}
With a closer analysis of the second term in ({\ref{TERMS}}) we will prove the following.  
\begin{proposition}[Decay of the second term in ({\ref{TERMS}}):  Closer analysis]
\label{otherdecay}
For all $r > 0$ and for the unique integer $K=K(r)$ such that $r_{K+1} < r \leq r_{K+2}$, we have
\[  \P \big( \obracket{E_K^c \cup C_K(r^ \a)} \cap \naiveset{ Z > r} \big) \leq {\constantmethoda}r^{\frac{-\a d}{2}} + \constantmethodb r^{ 2(\a-1)} (\log r)^{\oldd },  \]

for some $\constantmethoda = \constantmethoda(d) > 0$ and $\constantmethodb = \constantmethodb(d) > 0$.
\end{proposition}
\noindent Proposition {\ref{otherdecay}} together with Proposition {\ref{CLTdecay}} yields a proof of Theorem {\ref{result}}.  
\begin{proof}[Proof of Theorem {\ref{result}}]
From the previous results, ({\ref{TERMS}}), Proposition {\ref{CLTdecay}}, and Proposition {\ref{otherdecay}}  we have that
\[ \PPp{Z_{\Phi} > r} \leq \cclt r^{\frac{-\a d}{2}} + {\constantmethoda} r^{\frac{-\a d}{2}} + \constantmethodb r^{ 2(\a-1)} (\log r)^{\oldd }.\]
Hence we are led to minimize the quantity: $\max ( \frac{-\a d}{2},2(\a -1))$.
It is easy to verify that we should take $\a(d) = \choiceofalpha$.   Let $\beta(d):= d\a(d)/2 = \choiceofbeta$.  Thus we obtain that, 
\[ \result, \]
where $c = c(d) < \infty$ and $\conditions$.   
\end{proof}
We will now work towards a proof of Proposition {\ref{otherdecay}}.  We will need to examine the geometry of the blobs a bit closer to prove Proposition {\ref{otherdecay}}.  Again in light of ({\ref{encdecay}}) we do not need to worry about the event $E_k^c$.     Let us consider the decomposition,
\begin{equation}
\label{decomposition}
{  C_k(r^\a) \ = \ E_k^c \cap C_k(r ^ \a) \hspace{0.2cm}  \cup \hspace{0.2 cm}   E_k  \cap C_k(r ^ \a). }
\end{equation}
The second term puts us in a position akin to the situation of Proposition {\ref{CLTdecay}},  since we can control the diameter of the \kblob \ containing the origin.  We now examine two situations.  One where the \kblob \ containing $\Oo$ is possibly very small (see Lemma {\ref{ckb}} and Figure {\ref{smallblob}}) and another where there are {\em{enough}} points inside the \kblob \  to make good use of the central limit theorem \ (see Lemma {\ref{cka}} and Figure {\ref{bigblob}}).

Let $j \geq k$, consider again \jseeds \ on the sets: $$A_{j} = A_j(r^\a):=  S(-s_j, r_j + r^\a) \setminus S(-s_j, r_j - r^\a).$$  Observe that  seeds on two levels will not overlap; that is a seed in $A_{j}$ will not overlap with a seed in $A_{m}$, for $j \not = m$.  Also recall that by our definition of $k$-seeds, no two (distinct) $k$-seeds will overlap.  Since $C_k(r^\a)$ is the event that for some $j \geq k$, the set  $A_j$ contains a $j$-seed.  We will further split up this event.  
Define:
\begin{eqnarray}
\label{decomp}
  C^1_{k}(r^\a)&:=&\{ \text{for all $j \geq k$, the set $A_j$ contains at most one $j$-seed and} \nonumber \\ &\hspace{0.1 cm}&  \text{there is an unique $j \geq k$ such that $A_j$ contains a $j$-seed} \}   \nonumber \\
  C^2_{k}(r^\a)&:=& C_k(r^\a) \setminus C^1_{k}(r^\a)
\end{eqnarray}
We will throw away the term $\naiveset{Z > r}$  when we bound $\P(C^2_{k}(r^\a) \cap E_k \cap \naiveset{Z > r}    )   $, but we will keep it when we bound $\P(C^1_{k}(r^\a)   \cap E_k \cap \naiveset{Z > r} )$.  
\begin{lemma}
\label{ckb}
For all $r > 0$ and for the unique integer $K=K(r)$ such that $r_{K+1} < r \leq r_{K+2}$, we have $\PPp{C^2_{K}(r^\a)} \leq  \constantmethodbb (r^{ \a-1} (\log r)^{\old })^2$,
for some $\constantmethodbb = \constantmethodbb(d) > 0$.
\end{lemma}
\begin{proof}
For all $ j \geq 1,$ let  $U_j = \ns{A_j \ \text{contains} \ a \  j\text{-seed}}.$  Thus from ({\ref{AKK}}), we have $U_j = U_j(r^\a)$.  Since for all $j \geq 1$, no two (distinct) $j$-seeds overlap we   have that
\begin{equation*}
\PP{\text{$A_j$ contains more than one $j$-seed}} \leq \P(U_j)^2.
\end{equation*}
Similarly, since seeds in $A_j$ and $A_m$ do not overlap for $j \not = m,$ we have  
$\P(U_j \cap U_m) =  \P(U_j)\P(U_m)$, for all $j \not = m$.
Since $$C^2_k (r ^ \a) \ \subset \ \left( \bigcup_{j > m \geq k}  U_j \cap U_m \right) \ \  \bigcup \ \ \left( \bigcup_{j \geq m } \ns{\text{$A_j$ contains more than one $j$-seed}}\right),$$ we have $$\P(C^2_k(r^\a)) \leq \sum_{j \geq k} \P(U_j) \sum_{m \geq k} \P(U_m).$$  By equations ({\ref{evident}}) and (\ref{kdbound}) it is easy to see that $\PPp{C^2_{K}(r^\a)} \leq  \constantmethodbb (r^{ \a-1} (\log r)^{\old })^2,$ for some $\constantmethodbb = \constantmethodbb(d) > 0$.
\end{proof}

Thus we have an improved term $r^{2(\a -1)}$.  For a realization of the event $C^2_{k}(r^\a)$ see Figure {\ref{smallblob}}.
\begin{figure}[htbp]
\begin{center}
\includegraphics[width=3.5cm]{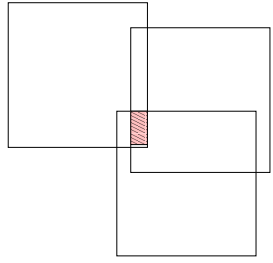}
\end{center}
\caption{The shaded region represents the $k$-blob containing the origin.  Notice that on the event $C ^{2}_k(r^\a)$ the $k$-blob can be quite {\em{small}}.  For this reason it seems we will not be able to do any better by including  $\naiveset{Z > r}$. }
\label{smallblob}
\end{figure}

We now turn our attention to the event $C^1_{k}(r^\a)$.
\begin{lemma}
\label{TERMCCC1}
For all $r > 0$ and for the unique integer $K=K(r)$ such that $r_{K+1} < r \leq r_{K+2}$, we have
$C^1_{K} (r^ \a)  \subset \naiveset{ \cardsupp{L_K} \geq \ccube{r^{\a d}}}$, for some constant, $0 < \ccube= \ccube(d) < \infty$.
\end{lemma}
\begin{proof}
On the event $C^1_{K}(r^\a)$, there is exactly one $j$-cutter that has the property that it intersects $S(r ^\a)$ and $j \geq K$; call this unique cutter $C$.  Observe that if the cutter $C$ was removed, the blob containing the origin would contain all of $S(r ^ \a)$.    Note that $C$ has side-length at least $2r_K$.   It is easy to see that there is a constant $\csidele < \infty$ independent of $k$ so that $\csidele r_k \geq r_{k+2}.$  Thus $\csidele r_K \geq r.$  Therefore the blob containing the origin must contain a $d$-cube with side-length $\frac{r^\a}{\csidele}.$  
\begin{figure}[htbp]
\begin{center}
\includegraphics[width=3.5cm]{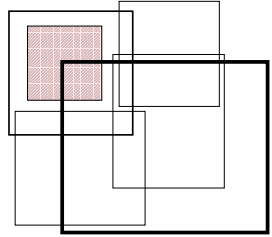}
\end{center}
\caption{An illustration of the event $C^1_{k}( r^\a) \cap E_k$.  The shaded region represents the restricted cutter zone of radius $r^\a$ about the origin.  The very thick cutter represents the unique $j \geq k$ that intersects $S(r^\a)$.  This corresponds to the event $C^1_{k}(r^\a)$.   The other thick cutter represents a $k$ cutter enclosing the origin.  This corresponds to the event $E_k$.  }
\label{bigblob}
\end{figure}
\end{proof}
\begin{lemma}
\label{cka}
For all $r > 0$ and for the unique integer $K=K(r)$ such that $r_{K+1} < r \leq r_{K+2}$, we have
 $\P \big( C^1_{K}(r^\a) \cap  E_K  \cap \naiveset{Z > r} \big)  \leq \frac{\constantmethodaa}{ (r^{\a})^ \frac{d}{2}   }$, for some $\constantmethodaa = \constantmethodaa(d) > 0$.
\end{lemma}
\begin{proof}
Again from ({\ref{TERME}}) and ({\ref{choiceofk}}) we have  $C ^1_{K}(r^\a) \cap  E_K  \cap \naiveset{Z > r} \subset \naiveset{ \Oo \textte{is \Kbad}}$.
So by Lemma {\ref{TERMCCC1}}, it  suffices to show that for all $k \geq 1$, we have
\begin{equation}
\label{condind2}
{
  \P \Big(   \naiveset{\Oo \ \text{is \kbad} }  \cap E_k  \cap \naiveset{ \cardsupp{L_k} \geq \ccube r^{\a d}   }    \Big) \leq \frac{\constantmethodaa}{(r^\a)^{\frac{d}{2}}},   
}
\end{equation}
from some $\constantmethodaa= \constantmethodaa(d) > 0$.
Equation \ref{condind2} follows from  Lemma {\ref{againuse}}.
\end{proof}
\begin{proof}[Proof of Proposition {\ref{otherdecay}}]
Using ({\ref{decomposition}}) and ({\ref{decomp}}), we have that
\begin{eqnarray*}
\P \big( \obracket{E_K^c \cup C_K(r^ \a)} \cap \naiveset{ Z > r} \big) &\leq&  
 \P(E_K ^c) + \P(C^2_K(r^\a)) + \\ 
&\hspace{0.1 cm}& \P \big( \obracket{C ^1_K(r^ \a)} \cap E_K \cap \naiveset{ Z > r} \big). 
\end{eqnarray*}
From equation (\ref{encdecay}) and  Lemmas {\ref{ckb}} and {\ref{cka}} we obtain that 
\begin{eqnarray*}
\P \big( \obracket{E_K^c \cup C_K(r^ \a)} \cap \naiveset{ Z > r} \big) &\leq&  
 \cestager r^{-\frac{d}{1 + \log 2}} + \constantmethodbb (r^{ \a-1} (\log r)^{\old })^2 + \\ 
&\hspace{0.1 cm}&  \frac{\constantmethodaa}{ (r^{\a})^ \frac{d}{2}   }. 
\end{eqnarray*}
\end{proof}

\section*{Open Problems}

\begin{enumerate}

\item
What is the optimal tail behaviour for  translation-equivariant matchings on $\Z^d$ in the case $d \geq 3$?  When $ d \geq 3$ from {\cite{Timarb}} for all $\e > 0$, there exists a translation-equivariant matching rule with exponential tails of order $\exp(-cr^{d-2-\e})$, where $c > 0$ is some constant.  Does there exists a translation-equivariant matching rule with tails of order $\exp(-cr^{d})$?  The original problem is from {\cite{Extra-Heads}} and it also contains a few other related open problems. 

\item
We say that a translation-equivariant matching rule is {\bf{oriented}} if it satisfies the additional restriction that if a site $x$ is matched to a site $y$ that contains a one,  then $y_i \geq x_i$ for all $i \leq d$.  Observe that in Meshalkin matching, a zero is always matched to a one that is to the right of it.   Note that it is not obvious that the method employed in this paper can be modified to work in an oriented setting.  In one dimension, the restriction of orientation does not make a difference; one might think it should not for higher dimensions as well.  What is the optimal tail behaviour for matchings in $\Z^d$ with the restriction that we consider orientation as well?

\end{enumerate}

\section*{Acknowledgements}
This research is funded in part by the NSERC and a U.B.C.\ Graduate Fellowship.  I like to thank my supervisor Ander Holroyd for introducing me to this topic.  Our numerous discussions and meetings were invaluable and his assistance and advice have helped to shape all aspects of this work.  I would also like to take the opportunity to thank Vlada Limic for introducing me to probability.

%\bibliography{references}
%\bibliographystyle{apt}

\end{document}